\documentclass{amsart}
\usepackage{amssymb}

\DeclareMathOperator{\adj}{adj}
\DeclareMathOperator{\Der}{Der}
\DeclareMathOperator{\id}{id}
\DeclareMathOperator{\tr}{tr}

\theoremstyle{definition}
\newtheorem{defn}{Definition}[section]

\newtheorem{remark}[defn]{Remark}
\newtheorem{sit}[defn]{}

\newtheorem{example}[defn]{Example}

\theoremstyle{plain}

\newtheorem{theorem}[defn]{Theorem}

\newtheorem{lemma}[defn]{Lemma}

\newtheorem{cor}[defn]{Corollary}
\newtheorem{corollary}[defn]{Corollary}

\theoremstyle{remark}

\begin{document}
\title[Factoring the Adjoint]{The Adjoint of an Even Size
Matrix Factors}

\author{Ragnar-Olaf Buchweitz}
\address{Dept.\ of Math., University of
Tor\-onto, Tor\-onto, Ont.\ M5S 3G3, Canada}
\email{ragnar@math.utoronto.ca}

\author{Graham Leuschke}
\address{Dept.\ of Math., University of
Tor\-onto, Tor\-onto, Ont.\ M5S 3G3, Canada}
\email{gleuschk@math.toronto.edu}

\thanks{The first author was partly supported by NSERC grant
3-642-114-80.}

\date{\today}

\begin{abstract}
We show that the adjoint matrix of a generic square matrix of 
even size can be factored nontrivially.  This answers a question of G.~Bergman.  This note should be considered a preliminary report on work in progress.
\end{abstract}
\maketitle

\section{Determinants and Derivations}
\begin{sit}
Let $K$ be a commutative ring, $X=(x_{ij})$ the generic
$(n\times n)$--matrix, whose entries thus form a family of
$n^{2}$ indeterminates, and set $S=K[x_{ij}]$, the
polynomial ring over $K$ in these variables.
\end{sit}

\begin{sit}
The determinant $\det(X)$ of the generic matrix $X$ is a
nonzerodivisor in $S$, and the classical adjoint matrix
$\adj(X)$ of $X$ is uniquely determined through either of
the following two matrix equations
\begin{equation}
\label{eq:1}
\adj(X)X = \det(X)\id_{n}\quad\text{and}\quad X\adj(X)= 
\det(X)\id_{n}\,,
\tag{$*$} 
\end{equation}
where $\id_{n}$ represents the $n\times n$ identity matrix.
\end{sit}

\begin{sit}
We will use the following notation for minors of the generic
matrix $X$: Let $[i_{1}i_{2}\ldots i_{k}\mid
j_{1}j_{2}\ldots j_{k}]$ denote the (unsigned) determinant
of the $(k\times k)$--submatrix of $X$ that consists of the
rows indexed $1\leq i_{1}<\cdots < i_{k}\leq n$, and of the
columns indexed $1\leq j_{1}<\cdots < j_{k}\leq n$.

The symbol $[i_{1}i_{2}\ldots i_{k}\;\widehat\mid\;
j_{1}j_{2}\ldots j_{k}]$ will denote the complementary
minor, thus, the determinant of the
$(n-k)\times(n-k)$--submatrix of $X$ obtained by removing
the rows indexed $i_{\nu}$ and the columns indexed
$j_{\nu}$.  For consistency, the empty determinant, for
$k=n$, has value $1$.

We extend the symbols $[?\mid\; ?]$ and $[?\;\widehat\mid\;?]$
to not necessarily strictly increasing index sets by
requiring them to be alternating in both the left and right
arguments. In particular, each symbol vanishes if there is 
repetition of indices either before or after the vertical 
bar.  
\end{sit}

\begin{sit} If $U$ is any $(n\times n)$--matrix over some
$K$--algebra $R$, then there exists a unique $K$--algebra
homomorphism $ev_{U}:S\to R, x_{ij}\mapsto u_{ij}$ that
transforms the entries of $X$ to those of $U$.  Let
$[\ldots](U) = ev_{U}([\ldots])$ represent the corresponding
minor of the matrix $U$, and write $I_{t}(U)\subseteq R$ for
the ideal generated by all the $(t\times t)$--minors of $U$.
The transpose of a matrix $U$ will be denoted $U^{T}$.
\end{sit}

\begin{example}
The $(i,j)$-th entry of the adjoint matrix can be written as
$$
\adj(X)_{ij} =
(-1)^{i+j}[j\;\widehat\mid\; i] = (-1)^{i+j}[1\ldots\widehat
j\ldots n\mid 1\ldots\widehat i\ldots n]\,.
$$
\end{example}

\begin{sit}
Recall that a map $D:R\to R$, on a not necessarily
commutative ring $R$, is a {\em derivation\/} if $D(ab) =
D(a)b + aD(b)$ for any elements $a,b\in R$.

For example, the {\em partial derivation\/} $\partial_{ij} =
\frac{\partial}{\partial x_{ij}}$ with respect to the
variable $x_{ij}$ defines a derivation on $S$ that is
furthermore $K$--linear.  These partial derivations form
indeed a basis of the free $S$--module $\Der_{K}(S)$ of all
$K$--linear derivations on $S$,
$$
\Der_{K}(S) \cong \bigoplus_{1\le i,j\le n}S\partial_{ij}\,.
$$
\end{sit}

Now we state the facts on derivations and minors that we 
will use.

\begin{lemma}
\label{lem:diffdet}
If $R$ is a commutative ring, $D:R\to R$ a derivation,
and $U$ an $(n\times n)$--matrix over $R$, then
$$
D(\det U) = \sum_{i=1}^{n}
\left|
\begin{matrix}
u_{11}&\cdots&u_{1n}\\
\vdots&&\vdots\\
D(u_{i1})&\cdots&D(u_{in})\\
\vdots&&\vdots\\
u_{n1}&\cdots&u_{nn}
\end{matrix}
\right|
= \sum_{i=1}^{n}
\left|
\begin{matrix}
u_{11}&\cdots&D(u_{1j})&\cdots&u_{1n}\\
\vdots&&\vdots&&\vdots\\
u_{n1}&\cdots&D(u_{nj})&\cdots&u_{nn}
\end{matrix}
\right|
$$
where $|V|$ denotes the determinant of the matrix $V$.
\end{lemma}
\begin{proof}
This follows immediately from the Leibnitz rule for
derivations applied to the complete expansion of the
determinant.
\end{proof}
\begin{lemma}
\label{lem:derdet} 
Let $X$ be again the generic matrix and $S$ the associated
polynomial ring over $K$.
\begin{enumerate}
\item[(1)]
\label{lem:derdet1}
For any pair of indices $1\le i,j\le n$,
$$
\partial_{ij}(\det X) = \adj(X)_{ji}
$$ 
equivalently,
$$
\adj(X)^{T} = (\partial_{ij}(\det X) )_{ij}\,.
$$

\item[(2)]
\label{lem:derdet2}
For any pair of indices $1\le i,j\le n$,
$$
\sum_{\nu=1}^{n}x_{i\nu}\partial_{j\nu}(\det X) =
\delta_{ij}\det(X) = \sum_{\nu=1}^{n}x_{\nu i}\partial_{\nu
j}(\det X)\,,
$$
where $\delta_{ij}$ is the {\em Kronecker symbol\/}.
\item[(3)]
\label{lem:derdet3}
For any indices $1\le i_{1}, i_{2},\ldots,i_{k}\le n$ and $1\le
j_{1},j_{2},\ldots,j_{k}\le n$, 
\begin{align*}
\partial_{i_{1}j_{1}}\cdots \partial_{i_{k}j_{k}}(\det X) &=
(-1)^{i_{1}+\cdots +i_{k}+j_{1}+\cdots j_{k}} [i_{1}\ldots
i_{k}\;\widehat\mid\, j_{1},\ldots,j_{k}]
\,,
\end{align*} 
in particular, these terms vanish whenever there is a 
repetition among the $i$'s or the $j$'s.
\end{enumerate}
\end{lemma}

\begin{proof}
Claim (1) follows from \ref{lem:diffdet} with
$D=\partial_{ij}$ and $U=X$.  In view of (1), claim (2) is
simply a reformulation of the equation (\ref{eq:1}) above.
To see (3), apply first \ref{lem:diffdet} or (1) to the generic
matrix using the derivation $\partial_{i_{k}j_{k}}$, and 
then use induction on $k\ge 1$.
\end{proof}
\section{The Factorizations}
We now use the ``differential calculus'' from the previous
section to establish two factorization results about
products of the transpose of the adjoint matrix with {\em
alternating\/} matrices on one or both sides.  Recall that
an $(n\times n)$--matrix $A=(a_{kl})$ is alternating if
$A^{T}= -A$ and the diagonal elements vanish, $a_{kk}=0$ for
each $k=1,\ldots, n$.  The latter condition is of course a
consequence of the first as soon as $2$ is a nonzerodivisor
in $K$.

\begin{theorem}
\label{thm:main}
Let $U, A$ be $(n\times n)$--matrices over a commutative
ring $R$, with $A$ alternating.  The $(n\times n)$--matrix
$B = (b_{rs})$ with entries from $I_{1}(A)\cdot I_{n-2}(U)
\subseteq R$, given by
$$
b_{rs}
= \sum_{k<l}a_{kl}(-1)^{k+l+r+s}[kl\;\widehat\mid\; rs](U)
$$
is then alternating as well and satisfies the matrix
equation
\begin{equation}
\label{eq:2}
\tag{$**$}
A\cdot \adj(U)^{T} = UB\,.
\end{equation}
If $\det U$ is a nonzerodivisor in $R$, then $B$ is the {\em
unique\/} solution to this equation.
\end{theorem}

\begin{proof}
As $[kl\;\widehat\mid\; sr] = -[kl\;\widehat\mid\; rs]$ and
$[kl\;\widehat\mid\; rr] = 0$, the matrix $B$ is 
alternating.  To verify that $B$ satisfies (\ref{eq:2}), it
suffices to establish the generic case, where $R=S$ and
$U=X$.  Let $E_{ij}$ denote the {\em elementary\/} $(n\times
n)$--matrix with $1$ at position $(i,j)$ as its only nonzero
entry.  Recall that $E_{rs}E_{uv}= \delta_{su}E_{rv}$ for
any indices $1\le r,s,u,v\le n$.  As
$\partial_{kr}\partial_{ls}(\det X) =
(-1)^{k+l+r+s}[kl\;\widehat\mid\; rs]$ by
Lemma~\ref{lem:derdet}(3), the right hand side of
(\ref{eq:2}) expands now first as
\begin{align*}
XB &= \bigg(\sum_{i,\nu}x_{i\nu}E_{i\nu}\bigg)
\bigg(\sum_{\mu,j}\sum_{k<l}a_{kl}\partial_{k\mu}
\partial_{lj}(\det X)E_{\mu j}\bigg)\\
&=\sum_{k<l}a_{kl}\sum_{i,j}\bigg(\sum_{\nu}x_{i\nu}
\partial_{k\nu}\partial_{lj}(\det X) \bigg)E_{ij}\,.
\end{align*}
The innermost sum can be simplified using first that partial
derivatives commute, 
then applying the product
rule, and finally invoking Lemma~\ref{lem:derdet}(2)
together with the fact that $\partial_{lj}(x_{i\nu}) =
\delta_{il}\delta_{j\nu}$. In detail, these steps yield the 
following equalities:
\begin{align*}
\sum_\nu x_{i\nu} \partial_{k\nu}\partial_{lj}(\det X)
&= \sum_\nu x_{i\nu} \partial_{lj}\partial_{k\nu}(\det X)\\
&= \sum_\nu \partial_{lj}\big(x_{i\nu} \partial_{k\nu}(\det
X)\big) - \sum_\nu
\partial_{lj}(x_{i\nu})\partial_{k\nu}(\det X) \\
&= \partial_{lj}\bigg(\sum_\nu x_{i\nu} \partial_{k\nu}(\det
X)\bigg) - \delta_{il} \sum_\nu \delta_{j\nu}
\partial_{k\nu}(\det X) \\
&= \delta_{ik} \partial_{lj} (\det X) -
  \delta_{il} \partial_{kj}(\det X) 
\end{align*}
In light of this simplification, we may expand $XB$ further 
as follows:
\begin{align*}
XB &=\sum_{k<l}a_{kl}\sum_{i,j}\bigg(\sum_{\nu}x_{i\nu}
\partial_{k\nu}\partial_{lj}(\det X) \bigg)E_{ij}\\
&= \sum_{k<l}a_{kl}\sum_{i,j}\bigg(\delta_{ik}
\partial_{lj}(\det X) - \delta_{il}\partial_{kj}(\det
X)\bigg)E_{ij}\\
&= \sum_{k<l}a_{kl}\sum_{j}\bigg(\partial_{lj}(\det X)E_{kj}
- \partial_{kj}(\det X)E_{lj}\bigg)\\
&= \sum_{k<l}a_{kl}\bigg(E_{kl}\sum_{j}\partial_{lj}(\det
X)E_{lj} -E_{lk} \sum_{j} \partial_{kj}(\det
X)E_{kj}\bigg)\\
&= \sum_{k<l}a_{kl}\bigg(E_{kl}\sum_{i,j}
\partial_{ij}(\det X)E_{ij} -E_{lk}
\sum_{i,j}\partial_{ij}(\det X)E_{ij}\bigg)\\
&= \sum_{k<l}a_{kl}\big(E_{kl}-E_{lk}\big)
\sum_{i,j}\partial_{ij}(\det X)E_{ij}\\
&= A\cdot \adj(X)^{T}
\end{align*}
with the last equality using that $A$ is alternating, thus
$A=\sum_{k<l}a_{kl}(E_{kl}-E_{lk})$, and that $\adj(X)^{T} =
\sum_{ij}\partial_{ij}(\det X)E_{ij}$, in view of
\ref{lem:derdet}(1).

The final assertion about uniqueness follows from
(\ref{eq:2}) by multiplying from the left with $\adj(U)$ and
using equation (\ref{eq:1}) to obtain
$$
\adj(U)\cdot A\cdot \adj(U)^{T} = \det(U) B\,.
$$
\end{proof}


\begin{remark} 
One may formulate \ref{thm:main} equally well for
multiplication of the transpose of the adjoint matrix from
the right by an alternating matrix.  Namely, assume $A,B$
are $(n\times n)$--matrices over $S$ satisfying
$A\adj(X)^{T} = XB$.  Let $\varphi :S\to S$ be the
$K$--algebra automorphism uniquely determined through
$\varphi(x_{ij})= x_{ji}$.  Clearly, $\varphi$ is involutive
and exchanges $X$ and its transpose, $\varphi(X) = X^{T}$.
Moreover, $\varphi(\adj(X)) = \adj(X)^{T}$, in view of
equation (\ref{eq:1}).  Now
\begin{alignat*}{2}
A\adj(X)^{T} &= XB&\quad&\text{if, and only if,}\\
\varphi(A)\varphi(\adj(X)^{T}) &=
\varphi(X)\varphi(B)&&\text{if, and only if,}\\
\varphi(A)\adj(X) &= X^{T}\varphi(B)&&\text{if, and only
if,}\\
\adj(X)^{T}\varphi(A)^{T} &= \varphi(B)^{T}X\,.
\end{alignat*}
In case $A, B$ are alternating, then so are $\varphi(A), 
\varphi(B)$ and the last equation is equivalent to
$$
\adj(X)^{T}\varphi(A) = \varphi(B)X\,.
$$
\end{remark}

  We now investigate 
what happens when multiplying simultaneously from both left
 and right.

\begin{theorem}
\label{thm:main2}
Let $U,A,B$ denote the same matrices as introduced in {\em
\ref{thm:main}}.  If $A'$ is another alternating $(n\times
n)$--matrix, then the $(n\times n)$--matrix $C=(c_{wm})$
with entries from $I_{1}(A)\cdot I_{n-3}(U)\cdot I_{1}(A')\subseteq R$
given by
\begin{align*}
c_{wm}&= \sum_{k<l, u<v}(-1)^{k+l+m+u+v+w}a_{kl}
[klm\;\widehat\mid\;uvw](U)a'_{uv}
\end{align*}
{satisfies}
\begin{equation}
\tag{***}
BA' = r\id_{n} + CU\,,
\end{equation}
{where}
\begin{align*}
r &= -\sum_{k<l, u<v}(-1)^{k+l+u+v}a_{kl}
[kl\;\widehat\mid\;uv](U)a'_{uv}\ \in R\,.
\end{align*}
\end{theorem}

\begin{proof}
It suffices again to verify the result for the generic 
matrix $U=X$, in which case we can employ once more the 
description of minors as given in \ref{lem:derdet}(3).
The straighforward calculation proceeds then as follows:
\begin{align*}
(BA' - r\id_{n})_{ij}&=
\sum_{m}\sum_{k<l}(-1)^{k+l+i+m}a_{kl}[kl\;\widehat\mid\;
im]a'_{m j} \\
&\quad + \delta_{ij}\sum_{k<l}\sum_{u<v}(-1)^{k+l+u+v}a_{kl}
[kl\;\widehat\mid\;uv]a'_{uv}\\
&=\sum_{k<l}a_{kl}\left(\sum_{m}\partial_{ki} 
\partial_{lm}(\det X)a'_{m j} + \sum_{u<v}\partial_{ku} 
\partial_{lv}(\det X)a'_{uv}\delta_{ij}\right)\\
&=\sum_{k<l,m}a_{kl}\left(\partial_{ki}
\partial_{lm}(\det X)a'_{m j} + \sum_{u<v}\partial_{ku}
\partial_{lv}\left(\partial_{m i}(\det X)x_{m
j}\right)a'_{uv}\right)
\end{align*}
where we have used \ref{lem:derdet}(2) in the last step.
Using the product rule twice together with
$
\partial_{rs}(x_{mn})=\delta_{rm}\delta_{sn}\,,
$
we find next
\begin{align*}
\partial_{ku} \partial_{lv}\left(\partial_{m i}(\det
X)x_{m j}\right)&= \partial_{ku}
\partial_{m i}(\det X)\delta_{lm}\delta_{vj}
+\partial_{lv}
\partial_{m i}(\det X)\delta_{km}\delta_{uj}\\
&\quad +\partial_{ku}
\partial_{lv}\partial_{m i}(\det X)x_{m j}
\end{align*}
Substituting and evaluating the Kronecker symbols yields
\begin{align*}
(BA' - r\id_{n})_{ij}
&=\sum_{k<l,m}a_{kl}\left(\partial_{ki}
\partial_{lm}(\det X)a'_{m j} + \sum_{u<v}\partial_{ku}
\partial_{lv}\left(\partial_{m i}(\det X)x_{m
j}\right)a'_{uv}\right)\\
&=\sum_{k<l}a_{kl}\left(\sum_{m}\partial_{ki}
\partial_{lm}(\det X)a'_{m j} +
\sum_{u<j}\partial_{ku}
\partial_{li}(\det X)a'_{uj}\right.\\
&\qquad\qquad  \left. +
\sum_{j<v}\partial_{ki}
\partial_{lv}(\det X)a'_{jv}
+ \sum_{u<v,m}\partial_{ku}
\partial_{lv}\partial_{m i}(\det X)a'_{uv}x_{m j}
\right)
\end{align*}
{The terms involving only second order derivatives
of the determinant cancel. To see this, rename
summation indices, use that
$
\partial_{km} \partial_{li}(\det
X)=-\partial_{ki} \partial_{lm}(\det X)
$ 
and that $A'$ is alternating, whence its entries satisfy
$a'_{mm}=0, a'_{jm}=-a'_{mj}$.  In detail,} 
\begin{align*}
(BA' - r\id_{n})_{ij}
&=\sum_{k<l}a_{kl}\left(\sum_{m}\partial_{ki}
\partial_{lm}(\det X)a'_{m j} - \sum_{m<j}\partial_{km}
\partial_{li}(\det X)a'_{m j}\right.\\
&\quad\quad  \left. -
\sum_{j<m}\partial_{ki}
\partial_{lm}(\det X)a'_{m j}
+ \sum_{u<v}\sum_{m}\partial_{ku}
\partial_{lv}\partial_{m i}(\det X)a'_{uv}x_{m j}
\right)\\
&=\sum_{m}\left(\sum_{
k<l,
u<v
}a_{kl} \partial_{ku}
\partial_{lv}\partial_{m i}(\det X)a'_{uv}\right)x_{m 
j}\\
&=\sum_{m}\left(\sum_{ k<l, u<v
}(-1)^{k+l+m+u+v+i}a_{kl}
[klm\;\widehat\mid\;uvi]a'_{uv}\right)x_{m j}\\
&= (CX)_{ij}
\end{align*}
where we have evaluated the third order derivatives of the
determinant according to \ref{lem:derdet}(3).
\end{proof}

Combining the results from  \ref{thm:main}
and \ref{thm:main2} yields the following.

\begin{corollary}
\label{cor:bifac}
Let $U,A,A'$ be $(n\times n)$--matrices over a commutative
ring $R$, with $A,A'$ {\em alternating}. One then has an 
equality of matrices
$$
A\adj(U)^{T}A' = rU + UCU\,,
$$
where $r$ and $C$ are as specified in {\em \ref{thm:main2}}.
\qed
\end{corollary}

\begin{remark}
The element $r\in I_{1}(A)\cdot I_{n-2}(U)\cdot I_{1}(A')\subseteq R$ is
a ``{\em half trace\/}'' of $BA'$, as
\begin{align*}
\tr(BA') &= \sum_{k<l}\sum_{i,j}a_{kl}(-1)^{k+l+i+j}[kl\;
\widehat\mid\;ij]a'_{ji}\\
&= 2\sum_{k<l}\sum_{i<j}a_{kl}(-1)^{k+l+i+j}[kl\;
\widehat\mid\;ij]a'_{ji}\\
&=2r
\end{align*}
invoking once again that $A'$ is alternating.  Equivalently, 
$\tr(CU) = (2-n)r\,.$
\end{remark}
\medskip

\begin{remark}
If $n=2$, all expressions of the form
$[klm\;\widehat\mid\;uvi]$ vanish, and \ref{thm:main}
together with \ref{thm:main2} specialize to the easily
established identity
\begin{align*}
\begin{pmatrix}
0&a\\
-a&0
\end{pmatrix}
\begin{pmatrix}
x_{22}&-x_{21}\\
-x_{12}&x_{11}
\end{pmatrix}
\begin{pmatrix}
0&b\\
-b&0
\end{pmatrix}
=
-ab
\begin{pmatrix}
x_{11}&x_{12}\\
x_{21}&x_{22}
\end{pmatrix}\,.
\end{align*}
\end{remark}
\medskip

If the size $n=2m$ is even, then there are {\em
invertible\/} alternating matrices of that size over any
commutative ring.  For example, the alternating ``hyperbolic
matrix'' 
$
\begin{pmatrix}
0&\id_{m}\\
-\id_{m}&0
\end{pmatrix}
$ has determinant equal to 1 over any ring.

\begin{cor} If $n$ is {\em even\/}, then the adjoint of the 
generic matrix admits nontrivial factorizations
$$
\adj(X) = YZ = Y'Z'
$$
into products of $(n\times n)$--matrices over 
$S$ with $\det(Y) = \det(Z') = \det(X)$.

More precisely, any pair of alternating $(n\times
n)$--matrices $A,A'$ of determinant equal to 1 over $S$
gives rise to such factorizations.  With $r$ and $C$ the
data associated to $A, A'$ as in {\em \ref{thm:main2}\/},
one may take
\begin{align*}
Y &= (A')^{-1}X^{T}\quad\text{and}\quad Z = 
(r\id_{n} + C^{T}X^{T})A^{-1}\,,\\
Y' &= (A')^{-1}(r\id_{n} +
X^{T}C^{T})\quad\text{and}\quad Z' = X^{T}A^{-1}\,.
\end{align*}
\end{cor}

\begin{proof}
Transposing the equation in \ref{cor:bifac} for $U=X$ yields
first
$$
(A')^{T} \adj(X) A^{T} = rX^{T} + X^{T}C^{T}X^{T}\,.
$$
As $A, A'$ are invertible and alternating, this equality is 
equivalent to
$$
\adj(X) = (A')^{-1}(rX^{T} + X^{T}C^{T}X^{T})A^{-1}\,.
$$
\end{proof}

\begin{remark}
Bergman \cite{Bergman:2003} shows that, over a field
$K$ of characteristic zero, in any factorization
$\adj(X) = YZ$ of the generic adjoint matrix into 
noninvertible factors, either $\det(Y) = \det(X)$ or
$\det(Z) = \det(X)$, up to units of $S\,$.
\end{remark}

\end{document}